\documentclass[leqno,a4paper]{article}

\usepackage{latexsym}
\usepackage{amssymb}
\usepackage{amsthm}
\usepackage{amsmath}

\theoremstyle{definition}

\newcommand{\tP}{\ensuremath{{}^tP}}

\renewcommand{\L}{\ensuremath{\mathrm{L}}}
\newcommand{\Con}{\ensuremath{\mathrm{C}}}
\newcommand{\Cinf}{\ensuremath{\mathrm{C}^\infty}}
\newcommand{\D}{\ensuremath{{\cal D}}}
\renewcommand{\S}{\ensuremath{{\cal S}}}
\newcommand{\E}{\ensuremath{{\cal E}}}

\newcommand{\mb}[1]{\ensuremath{\mathbb{#1}}}
\newcommand{\N}{\mb{N}}

\newcommand{\R}{\mb{R}}
\newcommand{\C}{\mb{C}}

\newcommand{\cl}[1]{\ensuremath{[#1]}}

\newcommand{\G}{\ensuremath{{\cal G}}}

\newcommand{\Gc}{\ensuremath{{\cal G}_\mathrm{c}}}
\newcommand{\EM}{\ensuremath{{\cal E}_{\mathrm{M}}}}

\newcommand{\NN}{\ensuremath{{\cal N}}}
\newcommand{\Ginf}{\ensuremath{\G^\infty}}

\newcommand{\WF}{\mathrm{WF}}
\newcommand{\singsupp}{\mathrm{sing supp}}
\newcommand{\Char}{\ensuremath{\text{Char}}}
\newcommand{\zs}{\setminus 0}
\newcommand{\CO}{\ensuremath{T^*(\Om) \zs}}

\renewcommand{\d}{\ensuremath{\partial}}

\newfont{\bl}{msbm10 scaled \magstep2}

\newtheorem{thm}{Theorem}[section]
\newtheorem{lemma}[thm]{Lemma}
\newtheorem{prop}[thm]{Proposition}

\newtheorem{cor}[thm]{Corollary}
\newtheorem{rem}[thm]{Remark}
\newtheorem{ex}[thm]{Example}
\newtheorem{asmpt}{Assumption}

\newcommand{\beq}{\begin{equation}}
\newcommand{\eeq}{\end{equation}}

\newcommand{\FT}[1]{\widehat{#1}}

\newcommand{\F}{\ensuremath{{\cal F}}}

\newcommand{\dis}[2]{\langle #1 , #2 \rangle}
\newcommand{\notmid}{\mid\kern-0.5em\not\kern0.5em}

\newcommand{\norm}[2]{{\| #1 \|}_{#2}}
\newcommand{\lone}[1]{\norm{#1}{\L^1}}

\newcommand{\linf}[1]{\norm{#1}{\L^\infty}}

\newcommand{\al}{\alpha}
\newcommand{\be}{\beta}
\newcommand{\ga}{\gamma}
\newcommand{\Ga}{\Gamma}
\newcommand{\de}{\delta}
\newcommand{\eps}{\varepsilon}

\newcommand{\vphi}{\varphi}

\newcommand{\la}{\lambda}
\newcommand{\La}{\Lambda}

\newcommand{\Om}{\Omega}
\newcommand{\sig}{\sigma}

\newcommand{\supp}{\mathop{\mathrm{supp}}}

\begin{document}

\title{Microlocal hypoellipticity 
of linear partial differential
operators with generalized functions as coefficients}

\author{\emph{G\"{u}nther H\"{o}rmann}\footnote{Supported by FWF grant P14576-MAT,
permanent affiliation: Institut f\"ur Mathematik, Universit\"at Wien}\\
\emph{Michael Oberguggenberger}\\
Institut f\"ur Technische Mathematik,\\ Geometrie und
Bauinformatik\\ Universit\"at Innsbruck\\
\ \\
\emph{Stevan Pilipovi\'{c}\ }\footnote{Supported by the MNTR of Serbia, 
Project 1835}\\
Institute of Mathematics and Informatics\\
Faculty of Science and Mathematics\\ University of Novi Sad
}

\date{}

\maketitle

\begin{abstract}
We investigate microlocal properties of partial differential operators with
generalized functions as coefficients. The main result is an extension of a
corresponding (microlocalized) distribution theoretic result on operators with
smooth hypoelliptic symbols. Methodological novelties and technical
refinements appear embedded into classical strategies of proof in order to
cope with most delicate interferences by non-smooth lower order terms. We
include simplified conditions which are applicable in special cases of
interest.
\end{abstract}


\section{Introduction}

Consider the linear partial differential equation
\[
    P u = f
\]
where $u$, $f$, as well as the coefficients of $P$ are generalized functions. We
investigate the general question of how to deduce information on the microlocal
regularity of $u$ from knowledge of the wave front set of $f$ and properties (of the
symbol) of $P$. As a matter of fact, a central issue in this is the identification of
appropriate conditions on $P$ which are generally applicable and, at the same time,
allow for nontrivial conclusions. If the coefficients of $P$ are smooth functions (or
analytic) and $u$ and $f$ are distributions (or ultradistributions or hyperfunctions)
then we may draw very detailed information revealed by the many efficient methods
from symbolic calculus, functional analysis, and microlocal analysis (cf.\
\cite{Hoermander:V1-4,Komatsu:91,Rodino:93,SKK:73,Taylor:81}).
  We take the freedom to refer to this, by now
well-established, distributional setting as the `classical case'. The
present paper is an exploration into a realm beyond, namely in the
direction where generalized functions may appear among the coefficients of
the operator $P$ (for example, to model jump discontinuities or even more
complex irregular properties of physical parameters). Since we are
concerned here with generally applicable information this puts us into facing
multiplication of distributions and thus obliges us to place our
analysis in the framework of algebras of generalized functions,
specifically Colombeau-type theories. A prominent new challenge in this
extended context is a more intricate mechanism to keep control over the
influence of non-smooth lower order terms. This is an observation stated
explicitly in \cite{HO:04} where results addressing this issue in special
cases were obtained. It is the purpose of this paper
to extend these results in substantial ways.

Subsection 1.1 serves as a review of the classical results which are to be extended
in the current paper. In Subsection 1.2 we give a brief account of the basics of the
modern theory of generalized functions in the sense of Colombeau. In particular,
notions of regularity in this setting are discussed in detail, including also a new
result on regular Colombeau functions with slow scale bounds. Section 2 constructs
the ``basic technical layer'' of the main proof in the paper. It provides also an
interface with ``higher levels'' of analytic conditions on the operators in terms of
two key assumptions. Section 3 presents the main result, the proof of which utilizes
the deductive structure laid out before. In Section 4 we give three types of
simplified conditions, which yield regularity results for certain special classes of
operators with non-smooth symbols. One of these is illustrated in the example of an
acoustic wave equation with bounded, but otherwise possibly highly irregular (e.g.\
discontinuous), coefficients.


\subsection{General regularity results for partial differential
operators with smooth coefficients}


To set the stage, we recall known, ``classical'' results relevant to the subsequent
generalizations. Let $\Om\subseteq\R^n$ be open and $\CO := \Om \times (\R^n\setminus
\{0\})$ be the cotangent space over $\Om$ with the zero section removed. Consider a
partial differential operator (PDO) $P$ of order $m$ with smooth coefficients,
given by $P(x,D) = \sum_{|\al| \leq m} a_\al(x) D^\al$, where
$a_\al\in\Cinf(\Om)$, and with principal part $P_m(x,D) = \sum_{|\al| = m}
a_\al(x) D^\al$. The (full) symbol of $P$ is $P(x,\xi) = \sum_{|\al| \leq m}
a_\al(x) \xi^\al$ and its principal symbol $P_m(x,\xi)$ is the symbol of
$P_m$. Both are interpreted as smooth functions on $T^*(\Om)$.

The notion of \emph{wave front set} of a distribution $u\in\D'(\Om)$ was introduced
in \cite{Hoermander:71}. We recall the elementary definition (due to
\cite[Proposition 2.5.5]{Hoermander:71}; see also \cite[Chapter 8]{Hoermander:V1-4})
of the wave front set of $\WF(u)\subseteq\CO$: The distribution $u$ is microlocally
regular at $(x,\xi)\in\CO$, i.e., $(x,\xi)\not\in\WF(u)$, if there exists a function
$\vphi\in\D(\Om)$ with $\vphi(x) \not= 0$ such that the Fourier transform of $\vphi
u$ is rapidly decreasing in a conic neighborhood of the direction $\xi$.

\paragraph{Noncharacteristic regularity:} Recall that actions of PDOs (with
smooth coefficients) on distributions do not enlarge singular supports. More
precisely, PDOs are \emph{microlocal} operators, which means that for any
$u\in\D'(\Om)$ we have \beq
    \WF(P u) \subseteq \WF(u).
\eeq The maximum failure of the reverse inclusion (in the general case) can be
captured by the ``geometry of the leading order terms'' of the operator. Recall that
$\Char(P) := P_m^{-1}(0) = \{ (x,\xi)\in \CO \mid P_m(x,\xi)=0 \}$ is the
characteristic set of $P$. It is a conic (with respect to the cotangent variable)
closed subset of $\CO$ and, together with $\WF(Pu)$, gives a general upper bound on
$\WF(u)$. This is the content of H\"ormander's theorem on noncharacteristic
regularity (cf.\ \cite[Theorem 8.3.1]{Hoermander:V1-4}).
\begin{thm} For any
$u\in\D'(\Om)$ we have the inclusion relation \beq \label{char_reg}
    \WF(u) \subseteq \WF(P u) \cup \Char(P).
\eeq
\end{thm}
The operator $P(x,D)$ is \emph{elliptic} if $\Char(P) = \emptyset$,
and \emph{(microlocally) elliptic at}
$(x_0,\xi_0)\in\CO$ if $P_m(x,\xi) \not=0$. Recall that microlocal ellipticity
at $(x_0,\xi_0)$ can be stated equivalently in terms of estimates on the full
symbol:
$\exists$ open set $U \ni x_0$ in $\Om$, $\exists$ open cone
    $\Ga \ni \xi_0$ in $\R^n\zs$, $\exists R > 0$, $C_0 > 0$
    $\forall \al, \be \in \N_0^n$, $|\be| \leq m$ $\exists C_{\al \be} > 0$
    such that
    \begin{align}
        |P(x,\xi)| & \geq C_0\, (1 + |\xi|)^m \label{ell_symb_a}\\
        |\d_x^\al \d_\xi^\be P(x,\xi)| & \leq C_{\al \be} \,
            |P(x,\xi)| \, (1 + |\xi|)^{-|\be|} \label{ell_symb_b}
    \end{align}
    for all $(x,\xi)\in U\times\Ga$ and $|\xi| \geq R$.

\begin{rem} 
(i) Recall that an operator $P(x,D)$ is said to be \emph{microhypoelliptic} if \[
\WF(u) = \WF(P u)\] for any distribution $u$, whereas it is \emph{hypoelliptic} if \[
\singsupp(u) = \singsupp(P u).\] Clearly, the former property implies the latter. If
$P(D)$ is a constant coefficient operator the two notions are equivalent (cf.\
\cite[Theorem 11.1.1]{Hoermander:V1-4}). But already operators with polynomial
coefficients can provide examples of differential operators which are hypoelliptic
and not microhypoelliptic (cf.\ \cite{PR:80}).  However, it follows from
(\ref{char_reg}) that all elliptic operators are microhypoelliptic. More generally,
pseudodifferential operators with hypoelliptic symbols are microhypoelliptic,
\cite[Chapter 22]{Hoermander:V1-4}, a result which we will state below in detail for
differential operators.

(ii) Any non-elliptic but microhypoelliptic operator shows that the
upper bound in (\ref{char_reg}) may be rather coarse. A simple example
is the heat operator with symbol $P(\xi,\tau) = i \tau + |\xi|^2$ with
$\Char(P) = \Om \times \{ (0,\tau) \mid \tau \not= 0 \}$.

\end{rem}

\paragraph{(Microlocal) Microhypoellipticity:}
Let $P(x,D)$ be a partial differential operator with smooth coefficients. It is said
to have a \emph{hypoelliptic symbol} $P(x,\xi)$ if the conditions
(\ref{hypo_symb_a}-\ref{hypo_symb_b}) below are satisfied on all of $\La = \CO$
(\cite[Definition 22.1.1]{Hoermander:V1-4}).

\begin{thm}
Let $(x_0,\xi_0)\in\CO$, $U\ni x_0$ open, $\Ga\ni\xi_0$ open conic (in
$\R^n\zs$), $m_0\in\R$, $0 \leq \delta < \rho \leq 1$, with the following
property: $\forall K \Subset U$ $\exists R > 0$ $\exists C_0 > 0$
$\forall \al, \be \in \N_0^n$, $|\be| \leq m$ $\exists C_{\al \be} > 0$
such that
\begin{align}
    |P(x,\xi)| & \geq C_0 \, (1+|\xi|)^{m_0} \label{hypo_symb_a}\\
    |\d_x^\al \d_\xi^\be P(x,\xi)| & \leq C_{\al \be}\,
        |P(x,\xi)|\, (1+|\xi|)^{-\rho |\be| + \de |\al|} \label{hypo_symb_b}
\end{align}
for all $(x,\xi)\in K\times \Ga$ with $|\xi| > R$. Then we have, with $\La =
U\times \Ga$, for any $u\in\D'(\Om)$
\beq \label{mhypo_reg_la}
    \La \cap \WF(u) = \La \cap \WF(P u).
\eeq
\end{thm}
(Cf.\ \cite[Theorem 22.1.4]{Hoermander:V1-4} for the global microhypoellipticity, and
\cite[Theorem 3.3.6]{MR:97} for the microlocal statement in the above form.)

\begin{cor} Let $M(P)$ be the union of all open conic subsets
$\La \subseteq \CO$ where conditions (\ref{hypo_symb_a}-\ref{hypo_symb_b}) are
satisfied ($m_0$, $\rho$, and $\de$ may vary with $\La$) then we have for any
distribution $u$ on $\Om$ the inclusion
\beq \label{mhypo_reg}
    \WF(u) \subseteq \WF(P u) \cup M(P)^c
\eeq
(the set theoretic complement taken in $\CO$).
\end{cor}

\begin{rem}
 If $M(P) = \CO$, i.e., $P$ has a hypoelliptic symbol,  then $P$ is
    microhypoelliptic, but the converse does not hold (cf.\
    \cite[Section 22.2]{Hoermander:V1-4} on generalized
    Kolmogorov equations). This is in contrast to the constant coefficient
    case where hypoelliptic symbols are in one-to-one correspondence with
    hypoelliptic operators (cf.\ \cite[Chapter 11]{Hoermander:V1-4}). Furthermore, in the
    latter case the set $M(P)$ can be determined by simple
    algebraic-geometric conditions (cf.\ \cite[p.15]{Hoermander:79a}).
In general, we have $M(P)^c \subseteq \Char(P)$ since
 (\ref{ell_symb_a}-\ref{ell_symb_b}) imply
(\ref{hypo_symb_a}-\ref{hypo_symb_b}). The inclusion can be strict, as the
    example of the heat operator, with $M(P)^c = \emptyset$
    shows. Therefore, (\ref{mhypo_reg}) is a refinement of
    (\ref{char_reg}).
\end{rem}

\subsection{Partial differential operators on algebras of
generalized functions}

\newcommand{\gR}{\widetilde{\mb{R}}}
\newcommand{\gC}{\widetilde{\mb{C}}}
\newcommand{\gOm}{\widetilde{\Omega}}
\newcommand{\gOmc}{\widetilde{\Omega}_\mathrm{c}}
\newcommand{\cO}{\ensuremath{{\cal O}}}

\subsubsection{Colombeau algebras}

The paper is placed in the framework of algebras of generalized
functions introduced by Colombeau in \cite{Colombeau:84,
Colombeau:85}.  We shall fix the notation and discuss a number of
known as well as new properties pertinent to Colombeau generalized
functions here. As a general reference we recommend \cite{GKOS:01}.

Let $\Omega$ be an open subset of $\R^n$. The basic objects of the
theory as we use it are families $(u_\eps)_{\eps \in (0,1]}$ of smooth
functions $u_\eps \in \Cinf(\Omega)$ for $0 < \eps \leq 1$. To
simplify the notation, we shall write $(u_\eps)_\eps$ in place of
$(u_\eps)_{\eps \in (0,1]}$ throughout.  We single out the following
subalgebras:
\vspace{2mm}\\
{\em Moderate families}, denoted by $\EM(\Omega)$, are defined by the property:
\begin{equation}
  \forall K \Subset \Omega\,\forall \alpha \in \N_0^n\,
    \exists p \geq 0:\;\sup_{x\in K} |\d^\alpha u_\eps(x)|
       = \cO(\eps^{-p})\ \rm{as}\ \eps \to 0\,.
     \label{mofu}
\end{equation}
{\em Null families}, denoted by $\NN(\Omega)$, are defined by the property:
\begin{equation}
  \forall K \Subset \Omega\,\forall \alpha \in \N_0^n\,
     \forall q \geq 0:\;\sup_{x\in K} |\d^\alpha u_\eps(x)|
        = \cO(\eps^q)\ \rm{as}\ \eps \to 0\,.
      \label{nufu}
\end{equation}
In words, moderate families satisfy a locally uniform polynomial estimate as $\eps
\to 0$, together with all derivatives, while null families vanish faster than any
power of $\eps$ in the same situation. The null families form a differential ideal in
the collection of moderate families. The {\em Colombeau algebra} is the factor
algebra
\[
   \G(\Omega) = \EM(\Omega)/\NN(\Omega)\,.
\]
The algebra $\G(\Omega)$ just defined coincides with the {\em special
Colombeau algebra} in \cite[Definition 1.2.2]{GKOS:01}, where the
notation $\G^s(\Omega)$ has been employed. However, as we will not use
other variants of the algebra, we drop the superscript $s$ in the
sequel.

Restrictions of the elements of $\G(\Omega)$ to open subsets of
$\Omega$ are defined on representatives in the obvious way. One can
show (see \cite[Theorem 1.2.4]{GKOS:01}) that $\Omega \to \G(\Omega)$
is a sheaf of differential algebras on $\R^n$. Thus the support of a
generalized function $u \in \G(\Omega)$ is well defined as the
complement of the largest open set on which $u$ vanishes. The
subalgebra of compactly supported Colombeau generalized functions will
be denoted by $\Gc(\Omega)$.

The space of compactly supported distributions is embedded in
$\G(\Omega)$ by convolution:
\[
   \iota:\E'(\Omega) \to \G(\Omega),\;
     \iota(w) = \cl{((w \ast \varphi_\eps)\vert_\Omega)_\eps} \,,
\]
where
\begin{equation}
   \varphi_\eps(x) = \eps^{-n}\varphi\left(x/\eps\right)
               \label{molli}
\end{equation}
is obtained by scaling a fixed test function $\varphi \in \S(\R^n)$ of integral one
with all moments vanishing. Here and henceforth the bracket notation $\cl{\ .\ }$ is
used to denote the equivalence class in $\G(\Om)$. By the sheaf property, this can be
extended in a unique way to an embedding of the space of distributions $\D'(\Omega)$,
and this embedding commutes with derivatives.

One of the main features of the Colombeau construction is the fact
that this embedding renders $\Cinf(\Omega)$ a faithful subalgebra.
(This property is optimal and cannot be improved to $\Con^k(\Om)$ for
finite $k$, in view of Schwartz' impossibility result
\cite{Schwartz:54}.) In fact, given $f \in \Cinf(\Omega)$, one can
define a corresponding element of $\G(\Omega)$ by the constant
embedding $\sigma(f) = \ \mbox{class of}\ [(\eps,x) \mapsto f(x)]$.
Then the important equality $\iota(f) = \sigma(f)$ holds in
$\G(\Omega)$. For a discussion of the overall properties of the
Colombeau algebra, we refer to the literature (e.g.
\cite{Colombeau:85, GKOS:01, O:92, Rosinger:90}).

Colombeau generalized numbers $\gC$ can be defined as the Colombeau algebra
$\G(\R^0)$, or alternatively as the ring of constants in $\G(\R^n)$. $\gC$ forms a
ring, but not a field. Concerning invertibility in $\gC$, we have the following
result (see \cite[Theorems 1.2.38 and 1.2.39]{GKOS:01}):

Let $r$ be an element of $\gR$ or $\gC$. Then
\begin{itemize}
\item[] $r$ is invertible if and only if
\item[] there exists some representative $(r_\eps)_\eps$
and an $m\in \N$ with $|r_\eps| \geq  \eps^m$ for sufficiently small
$\eps > 0$, if and only if
\item[] $r$ is not a zero divisor.
\end{itemize}

Concerning invertibility of Colombeau generalized functions, we may
state (see \cite[Theorem 1.2.5]{GKOS:01}):

Let $u \in \G(\Omega)$. Then
\begin{itemize}
\item[] $u$ possesses a multiplicative inverse if and only if
\item[] there exists some representative $(u_\eps)_\eps$ such that
for every compact set $K \subset \Omega$, there is $m\in \N$ with
$\inf_{x \in K}|u_\eps(x)| \geq  \eps^m$ for sufficiently small $\eps
> 0$.
\end{itemize}

In order to be able to speak about symbols of differential operators,
we shall need the notion of a polynomial with generalized
coefficients. The most straightforward definition is to consider a
generalized polynomial of degree $m$ as a member
\[
   \sum_{|\gamma| \leq m} a_\gamma \xi^\gamma \in \G_m[\xi]
\]
of the space of polynomials of degree $m$ in the indeterminate $\xi =
(\xi_1,\dots, \xi_n)$, with coefficients in $\G =
\G(\Omega)$. Alternatively, we can and will view $\G_m[\xi]$ as the
factor space
\begin{equation}\label{generalizedpolynomials}
    \G_m[\xi] = {\cal E}_{\mathrm{M},m}[\xi]/\NN_m[\xi]
\end{equation}
of families of polynomials of degree $m$ with moderate coefficients
modulo those with null coefficients. In this interpretation,
generalized polynomials $P(x,\xi)$ are represented by families
\[
   (P_\eps(x,\xi))_\eps = \Big(\sum_{|\gamma| \leq m} a_{\eps\gamma}(x) \xi^\gamma\Big)_\eps\,.
\]
Sometimes it will also be useful to regard polynomials as polynomial
functions and hence as elements of $\G(\Omega \times \R^n)$. Important
special cases are the polynomials with regular coefficients,
$\G_m^\infty[\xi]$, and with constant generalized coefficients,
$\gC_m[\xi]$. The union of the spaces of polynomials of degree $m$ are
the rings of polynomials $\G[\xi], \G^\infty[\xi]$, and
$\gC[\xi]$. Letting $D = (-i\d_1,\dots, -i\d_n)$, a differential
operator $P(x,D)$ with coefficients in $\G(\Omega)$ simply is an
element of $\G[D]$.

\subsubsection{Regularity of Colombeau functions}

We recall a few notions and present a new result about regularity of Colombeau
generalized functions.  In this setting, the notion is based on the subalgebra
$\G^\infty(\Omega)$ of {\em regular generalized functions} in $\G(\Omega)$. It is
defined by those elements which have a representative satisfying
\begin{eqnarray}
  \forall K \Subset \Omega\,\exists p \geq 0\,\forall \alpha \in \N_0^n:\;
  \sup_{x\in K} |\d^\alpha u_\eps(x)| = \cO(\eps^{-p})\
     \ \rm{as}\ \eps \to 0\,. \label{regefu}
\end{eqnarray}
Observe the change of quantifiers with respect to formula (\ref{mofu}); locally, all
derivatives of a regular generalized function have the same order of growth in $\eps
> 0$. One has that (see \cite[Theorem 25.2]{O:92})
\[
  \G^\infty(\Omega) \cap \D'(\Omega) = \Cinf(\Omega)\,.
\]
For the purpose of describing the regularity of Colombeau generalized
functions, $\G^\infty(\Omega)$ plays the same role as $\Cinf(\Omega)$
does in the setting of distributions.

The concept of \emph{microlocal regularity} of a Colombeau function follows the
classical idea of employing additional spectral information on the singularity from
the (Fourier) frequency domain (cf.\ \cite{Hoermann:99,HK:01,NPS:98}). It refines
$\Ginf$-regularity in the sense that the projection of the (generalized) wave front
set into the base space equals the (generalized) singular support. We recall the
definition of the generalized wave front set. First, $u\in\G(\Om)$ is said to be
\emph{microlocally regular} at $(x_0,\xi_0)\in\CO$ if (for a representative
$(u_\eps)_\eps\in\EM(\Om)$) there is an open neighborhood $U$ of $x_0$ and a conic
neighborhood $\Ga$ of $\xi_0$ such that for all $\vphi\in\D(U)$ we have that
$\F(\vphi u)$ is rapidly decreasing in $\Ga$; that means that $\exists N\in\R$
$\forall l\in\N_0$ $\exists C > 0$ $\exists \eps_0 > 0$:
\begin{equation}\label{rap_dec}
  |(\vphi u_\eps)\FT{\ }(\xi)| \leq C
  \eps^{-N} (1 + |\xi|)^{-l} \qquad \forall\xi\in\Ga,
        \forall\eps\in (0,\eps_0).
\end{equation}
(We denoted by $\FT{\ }$ the Fourier transform on test functions
and by $\F(\vphi u)$ the corresponding generalized Fourier transform
of the compactly supported Colombeau function $\vphi u$.) The
\emph{generalized wave front set} of $u$, denoted by $\WF_g(u)$, is defined as
the complement (in $\CO$) of the set of pairs $(x_0,\xi_0)$ where $u$ is
microlocally regular.

Let us recall a recently introduced notion which turns out to be crucial in the
context of regularity theory (cf.\ \cite{HO:04}), namely {\em slow scale nets}. By
this, we mean a moderate net of complex numbers $r = (r_\eps)_\eps \in \gC_M$ with
the following property: $\forall t \geq 0$, $\exists \eps_t > 0$ such that
$|r_\eps|^t \leq \eps^{-1}$ for all $\eps \in (0,\eps_t)$. Equivalently, defining the
{\em order of} $r$ by $\kappa(r) = \sup \{q \in \R: \exists \eps_q\,\exists C_q > 0:
|r_\eps| \leq C_q \eps^q\,, \forall \eps \in (0,\eps_q)\}$, $r$ is a slow scale net
if (and only if) it has order $\kappa(r) \geq 0$. We refer to \cite[Section 2]{HO:04}
for a detailed discussion of further properties.

An interesting property of the elements of $\G^\infty(\Omega)$ is that
their representatives are never bounded, unless all derivatives are
slow scale. More precisely, we have the following result; we denote by
$\EM^\infty(\Omega)$ the nets of smooth functions satisfying the
estimates (\ref{regefu}).
\begin{prop}\label{sc_ginf_prop}
Let $(u_\eps)_\eps \in \EM^\infty(\Omega), K \Subset \Omega$, and assume that
$\sup_{x \in K}|u_\eps(x)|$ is slow scale. Then $\sup_{x \in K}|\d^\alpha
u_\eps(x)|$ is slow scale for all $\alpha \in \N_0^n$.
\end{prop}
\begin{proof}
We can find a bounded open set $\omega \subset \Omega$ which is of cone
type and contains $K$; in fact, we can take $\omega$ as a finite union of
$n$-dimensional intervals with positive distance to the boundary of $\Omega$.
We use interpolation theory for the Sobolev spaces $H^\ell(\omega)$
to observe that $(k, \ell \geq 0)$
\[
    \big[L^2(\omega), H^{k + \ell}(\omega)\big]_{\ell/(k + \ell)} =
        H^\ell(\omega)\,,
\]
and for $v \in H^\ell(\omega)$,
\[
   ||v||_{H^\ell} \leq
        C_{k,\ell}||v||_{L^2}^{1 - \ell/(k + \ell)}
            ||v||_{k + \ell}^{\ell/(k + \ell)}\,,
\]
see e.~g. \cite[2.4.2, 4.3.1]{Triebel:78}. In particular,
\begin{equation} \label{interpolation}
   ||v||_{H^\ell}^{k + \ell} \leq
        C_{k,\ell}||v||_{L^2}^{k} ||v||_{k + \ell}^{\ell}\,.
\end{equation}
Now assume that $\sup_{x \in K}|\d^\alpha u_\eps(x)|$ is not slow scale for
some $\alpha, |\alpha| \geq 1$. Then its order is less than zero,
so there is $p > 0$ such that $\sup_{x \in K}|\d^\alpha
u_\eps(x)| > \eps^{-p}$ for a certain subsequence of $\eps \to 0$. By
Sobolev's embedding theorem, it follows that
\[
   ||u_\eps||_{H^\ell(\omega)} > C\eps^{-p}
\]
as well, provided $\ell > |\alpha| + n/2$, with some constant $C > 0$. On the
other hand,
\[
   ||u_\eps||_{L^2(\omega)} \leq
    C_\omega ||u_\eps||_{L^\infty(\omega)} =: r_\eps
\]
with $(r_\eps)_\eps$ slow scale. Inserting this in (\ref{interpolation}), we
get
\[
   \left(C\eps^{-p}\right)^{k + \ell} <
        C_{k,\ell}r_\eps^k||u_\eps||_{H^{k+\ell}(\omega)}^\ell\,,
    \forall k \geq 0\,.
\]
Thus
\[
   ||u_\eps||_{H^{k+\ell}(\omega)} >
        D_{k,\ell}\eps^{-p(1+k/\ell)}r_\eps^{-k/\ell},
\]
still for a subsequence as $\eps \to 0$. Given $k \geq 0$, we have that
$r_\eps^{-k/\ell}\eps^{-p} > 1$ for sufficiently small $\eps > 0$. Thus we get
that, given $k \geq =0$,
\[
   ||u_\eps||_{H^{k+\ell}(\omega)} > D_{k,\ell}\eps^{-pk/\ell)}
\]
for a subsequence of $\eps \to 0$. Since
\[
   ||u_\eps||_{W^{k+\ell,\infty}(\omega)} \geq
    C_\omega||u_\eps||_{H^{k+\ell}(\omega)}
\]
and $k$ is arbitrary, this contradicts $(u_\eps)_\eps \in \EM^\infty(\Omega)$
and proves the claim.
\end{proof}
\begin{cor}\label{sc_ginf_cor}
Let $u \in \G^\infty(\Omega), K \Subset \Omega$ and assume that for some
representative, $\sup_{x \in K}|u_\eps(x)|$ is bounded as $\eps \to 0$. Then
all derivatives of $u_\eps$ are slow scale on $K$.
\end{cor}
In particular, a locally bounded, regular generalized function has the property
that none of its derivatives can have a strictly negative order on any
compact set.
\begin{ex}
Polynomials of degree $\geq 1$ in the variable $x/\eps$ are elements of
$\EM(\Omega)$ which are not slow scale and attain negative orders in their
derivatives. On the other hand, typical examples of bounded generalized
functions are provided by regularizations of the Heaviside function: Let
$\varphi \in \S(\R)$ with integral one and define the mollifier
$\varphi_{r_\eps}$ as in (\ref{molli}). Then
$\big(H\ast\varphi_{r_\eps}\big)_\eps$ belongs to $\EM^\infty(\R)$ if and only
if $(1/r_\eps)_\eps$ is slow scale.
\end{ex}

The focus of the current work is on microlocalization of the notion of
hypoellipticity for PDOs with Colombeau functions as coefficients, as it was
introduced in \cite{HO:04}. Operators that enjoy the property
\[
   \big[u \in \G(\Omega)\,,\;  f \in \Ginf(\Omega) \mbox{ and }
     P(x,D)u =f \mbox{ in } \G(\Omega) \big]
    \Longrightarrow u \in \Ginf(\Omega)\,.
\]
 on every open subset
$\Omega \subset \R^n$ are called {\it $\Ginf$-hypoelliptic}. General results
on global elliptic regularity for operators with generalized constant
coefficients as well as on microlocal regularity for certain first-order
operators were obtained in \cite{HO:04}. It is also a source for a
variety of examples illustrating the above as well as related notions.
Recent related research in Colombeau regularity theory has shown progress in a
diversity of directions, including such topics as pseudodifferential
operators with generalized symbols and case studies in microlocal analysis of
nonlinear singularity propagation \cite{Garetto:04,HdH:01,HK:01,NPS:98}.

\section{The basic scheme: deducing regularity of the solution
from approximative solutions of the adjoint equation}

The proof of (\ref{char_reg}) in \cite[Theorem 8.3.1]{Hoermander:V1-4} is based on an
idea to use approximate solutions of the adjoint equation (with a particular
right-hand side) in deducing regularity of a solution to the original PDE. This
approach is elementary in the sense that it does not rely on pseudodifferential
operator machinery, though basic techniques naturally appear in it at an embryonic
stage. We adapt this procedure to $\G(\Om)$ in the constructions carried out below. A
closely related path was taken up already in \cite{DPS:98} and the current exposition
partly serves to correct and improve the results stated there.

\paragraph{Informal description of the underlying idea:} We briefly sketch the
strategy of the ``classical'' proof, i.e., in the context of operators with
smooth coefficients and distributional
solutions. Put $f = Pu$ and let $\vphi\in\D$, $\vphi(x_0) = 1$.
To deduce that $(x_0,\xi_0)\not\in \WF(u)$ one has to estimate
\beq \label{FT}
    (u \vphi)\FT{\ }(\xi) = \dis{u}{\vphi\, e^{- i\xi.}}
\eeq
for $\xi$ varying in a conic neighborhood of $x_0$. If the adjoint equation
\beq \label{adj_equ}
    \tP(\psi\, e^{-i\xi.}) = \vphi\, e^{-i\xi.}
\eeq
were solvable for some $\psi\in\D$ (with $\supp(\psi)$ close to
$\supp(\vphi)$)
then (\ref{FT}) could be rewritten as
\[
    \dis{u}{\tP(\psi e^{-i\xi.})} = \dis{f}{\psi e^{-i\xi.}} =
    (f \psi)\FT{\ }(\xi).
\]
Under the assumption $(x_0,\xi_0)\not\in \WF(f)$ this would prove the claim
and establish a relation of the type (\ref{char_reg}).

\paragraph{Skeleton of the microlocal regularity proof:} Let $u =
\cl{(u_\eps)_\eps}\in\G(\Om)$ and put $f_\eps = P_\eps u_\eps$, $f =
\cl{(f_\eps)_\eps}$. Assume $(x_0,\xi_0) \not\in \WF_g(f)$ and let
$U\times\Ga$ be an open conic neighborhood of $(x_0,\xi_0)$ such that for any
$\vphi\in\D(U)$ with $\vphi(x_0) = 1$ the Fourier transform $\F(\vphi f)$ is
rapidly decreasing in $\Ga$. This means that $\exists
M\in\N_0$ $\forall l \in\N_0$ $\exists C > 0$, $0 < \eps_0 < 1$ such
that
\beq\label{f_rap_dec}
    |(\vphi f_\eps)\FT{\ }(\xi)| \leq C \eps^{-M} (1 + |\xi|)^{-l}
    \qquad \forall \xi\in\Gamma, 0 < \eps < \eps_0.
\eeq

We want to show that
\beq\label{u_FT}
    (\vphi u_\eps)\FT{\ }(\xi) = \int u_\eps(x) \vphi(x) e^{-i\xi x}\, dx
\eeq
is rapidly decreasing in some conic neighborhood $\Ga_0 \subseteq \Ga$ of
$\xi_0$ when certain conditions on the family of symbols $P_\eps(x,\xi)$
($\eps\in (0,1]$) are met in $\supp(\vphi) \times \Ga$.

Note that it suffices to establish an estimate of the form (\ref{f_rap_dec})
whenever $|\xi| \geq r_\eps$, where $r_\eps >0$ is of slow scale and may
depend on $N$ (we refer to a corresponding remark in \cite[Section 6]{HO:04}).

We postpone the detailed discussion of various conditions on $P_\eps$
suiting the same proof skeleton until the following section. For the moment,
we will instead state general assumptions tailored directly towards the adjoint
operator method and investigate later on how they can be met in certain
circumstances.

The following lemma provides the basic algebraic mechanism in the construction
of approximate solutions of the adjoint equation. We re-investigate and state
the classical computations here in all details in order to prepare for
the close examination of the interplay of $\xi$-order and $\eps$-asymptotics
required later on.  One may think of the part of the operator $A$ in it to be played by $P_\eps$ (or $P_{\eps,m}$) and $Q$
corresponding to the adjoint $\tP_\eps$.
\begin{lemma}\label{alg_lem}
Let $A(x,D)$ and $Q(x,D)$ be partial differential operators of order $m$ on an
open subset $U$ of $\R^n$. Assume that  $A(x,\xi) \not= 0$ for all $(x,\xi)\in
U\times\Ga$. Then for any $w\in\Cinf(U)$ the following equation holds (on
$U\times\Ga$)
\beq\label{alg_equ}
    Q(x,D)\Big( \frac{e^{-i\xi x} w(x)}{A(x,\xi)} \Big) =
        e^{-i\xi x} \big( w(x) - R(\xi;x,D)w (x) \big)
\eeq
where $R(\xi;x,D) = - \sum_{|\be| \leq m} r_\be(x,\xi) D_x^\be$ is a partial
differential operator of order at most $m$ and coefficients given by
\begin{align}
    r_0(x,\xi) & = \frac{Q(x,-\xi) - A(x,\xi)}{A(x,\xi)} +
        \sum_{1 \leq |\ga| \leq m} \frac{1}{\ga !}\,
            \d_\xi^\ga Q(x,-\xi)\, D_x^\ga\big(\frac{1}{A(x,\xi)}\big)
        \label{r0}\\
    r_\be(x,\xi) & = \sum_{|\ga| \leq m-|\be|} \frac{1}{\be ! \ga !}\,
        \d_\xi^{\be + \ga} Q(x,-\xi)\, D_x^\ga\big(\frac{1}{A(x,\xi)}\big)
                \qquad  |\be| \geq 1. \label{rb}
\end{align}
\end{lemma}
\begin{proof} Using the fact that $\d_\xi^\be Q(x,D)(e^{-i\xi x}) =
e^{-i\xi x} \d_\xi^\be Q(x,-\xi)$ and by repeated application of the
H\"ormander-Leibniz formula (\cite[(1.1.10)]{Hoermander:V1-4}) we have
\begin{multline*}
    Q(x,D)\Big( \frac{e^{-i\xi x} w(x)}{A(x,\xi)} \Big) = \sum_{|\be|\leq m}
        \frac{1}{\be !} \,\d_\xi^\be Q(x,D)\big( \frac{e^{-i\xi x}}{A(x,\xi)}
        \big)\, D_x^\be w(x) \\
    = Q(x,D)\big(\frac{e^{-i\xi x}}{A(x,\xi)}
        \big)\, w(x) + \sum_{1 \leq |\be|\leq m} \frac{1}{\be !}
        \,\d_\xi^\be Q(x,D)\big( \frac{e^{-i\xi x}}{A(x,\xi)}
        \big)x \, D_x^\be w(x)\\
    = e^{-i\xi x} \Big(\sum_{|\ga| \leq m} \frac{1}{\ga !}\,
            \d_\xi^\ga Q(x,-\xi)\,
            D_x^\ga\big(\frac{1}{A(x,\xi)}\big) \, w(x) \\
    + \sum_{1 \leq |\be|\leq m} \sum_{|\ga| \leq m - |\be|}
        \frac{1}{\be ! \ga !}
        \,\d_\xi^{\be+\ga} Q(x,-\xi)\, D_x^\ga \big(\frac{1}{A(x,\xi)}\big)
        \, D_x^\be w(x)  \Big).
\end{multline*}
Here, (\ref{rb}) can be read off directly from the second (double) sum and
(\ref{r0}) follows by separating the term corresponding to $|\ga|=0$,
$Q(x,-\xi)/A(x,\xi)$, and adding and subtracting $A(x,\xi)w(x)/A(x,\xi)$.
\end{proof}

Formula (\ref{alg_equ}) is the basis for solving the adjoint equation
(\ref{adj_equ}) approximately.  For this, we apply the lemma to
$A(x,D) = P_\eps(x,D)$ and $Q(x,D) = \tP_\eps(x,D)$, where $\eps$ is
fixed but arbitrary in some interval $(0,\eps_0)$, $x\in U \subseteq \Om$
open, independent of    $\eps$, and
$\xi\in\Ga_\eps = \Ga \cap \{ |\eta| \geq r_\eps \}$, where the cone $\Ga$
is independent of $\eps$ and $r_\eps$ is of slow scale.

This defines a corresponding family $R_\eps(\xi;x,D)$ ($0 < \eps < \eps_0$,
$\xi\in\Ga_\eps$) of differential operators on $U$, with coefficients given by
corresponding parametrized versions of (\ref{r0}-\ref{rb}).  Observe that
$Q_m(x,-\xi) = P_{\eps,m}(x,\xi)$ and hence
$Q(x,-\xi) - A(x,\xi) = \tP_\eps(x,-\xi) - P_\eps(x,\xi)$ is a polynomial of
order at most $m-1$ with respect to $\xi$.

The second ingredient is the choice of $w$, which has to be linked with
$\vphi$ and will also depend on the parameters $\xi$ and $\eps$, as
well as on an approximation order $N\in\N$. We express this in the notation
$w_\eps^N(x,\xi)$. If we define
\beq\label{w_def}
    w^N_\eps(x,\xi) = \sum_{k=0}^{N-1} R_\eps(\xi;x,D)^k \vphi(x)
\eeq
then the expression $w_\eps^N - R_\eps w_\eps^N$ appearing on the right-hand
side of (\ref{alg_equ}) is a telescope sum and reduces to $\vphi - R_\eps^N
\vphi$. Therefore, equation (\ref{alg_equ}) yields in this case
\beq\label{appr_adj}
    \tP_\eps(x,D)\big(
        \frac{w_\eps^N(x,\xi)}{P_\eps(x,\xi)}e^{-i\xi x}\big)
    = e^{-i\xi x}\cdot \vphi(x) - e^{-i\xi x}\cdot R_\eps(\xi;x,D)^N \vphi(x)
\eeq
which is as close as we get to the informal requirement of (\ref{adj_equ}).

Equation (\ref{appr_adj}) suggests that
\beq\label{psi_def}
    \psi_\eps^N := w_\eps^N / P_\eps
\eeq
will give a
reasonable approximate solution in terms of decrease properties with
respect to $\xi$ (while keeping control over the $\eps$-growth), if the
operator family $(R_\eps)_\eps$ satisfies corresponding estimates. Note that
each $R_\eps(\xi;x,D)$ is a differential operator of order at most $m$ with
coefficients depending smoothly on $x\in U$ and being rational functions of
$\xi\in\Ga$.
\begin{asmpt}\label{R_asmpt} There is $M_1\in\N_0$ and $\tau > 0$ with the
property that $\forall N\in\N$ $\exists C > 0$,
$0 < \eps_1 \leq \eps_0$ and $\exists r_\eps >0$ of
slow scale such that
\beq\label{R_cond}
    |R_\eps(\xi;x,D)^N \vphi(x)| \leq C \eps^{-M_1} (1 + |\xi|)^{-N\tau}
\eeq
for all $x\in\supp(\vphi)$, $\xi\in\Ga$ with $|\xi| \geq r_\eps$, $\eps\in
(0,\eps_1)$.
\end{asmpt}

We deduce from (\ref{u_FT}) and
(\ref{appr_adj}) that
\begin{multline}\label{J_I_split}
     (\vphi u_\eps)\FT{\ }(\xi) =
    \int P_\eps(x,D)u_\eps(x) \cdot \psi_\eps^N(x,\xi)\, e^{-i\xi x}\, dx  \\
    + \int  u_\eps(x) e^{-i\xi x}\cdot R_\eps(\xi;x,D)^N \vphi(x)\, dx =:
    J_\eps^N(\xi) + I_\eps^N(\xi)
\end{multline}

\begin{lemma}\label{R_lemma}
If $(R_\eps)_\eps$ satisfies Assumption \ref{R_asmpt} then there exists
$N_1\in\N_0$ such that for all $N$ the integral $I_\eps^N(\xi)$
satisfies the following estimate:
$\exists C >0 $ $\exists \eps_2 > 0$
\beq\label{I_est}
   |I_\eps^N(\xi)| \leq C \eps^{- N_1}\, (1 + |\xi|)^{-N\tau}
\eeq
for all $x\in U$, $\xi\in\Ga$ with $|\xi| \geq r_\eps$, $0 < \eps < \eps_2$.
\end{lemma}
\begin{proof}
 Let $p\in\N_0$ and
$\eps_2\leq \eps_1$ be sufficiently small so that $|u_\eps(x)| \leq \eps^{-p}$
uniformly for $x\in\supp(\vphi)$ when $0 < \eps <\eps_2$. Then
\beq
   |I_\eps^N(\xi)| \leq C_\vphi\, \eps^{- M_1 - p}\, (1 + |\xi|)^{-N\tau}
\eeq
for all $x\in U$, $\xi\in\Ga$ with $|\xi| \geq r_\eps$, $0 < \eps < \eps_2$,
where $C_\vphi$ is the product of the constant in (\ref{R_cond}) and the
measure of $\supp(\vphi)$.
\end{proof}
Note that since $N_1$ is independent of $N\in\N$, (\ref{I_est}) can be
used in proving rapid decrease in (\ref{J_I_split})
once $J_\eps^N$ was shown to be rapidly decreasing. We observe that
\beq\label{J_int}
    J_\eps^N(\xi) = \int f_\eps(x)\,\psi_\eps^N(x,\xi)\, e^{-i\xi x}\, dx =
    \F_{x \to \nu}\big( f_\eps(x) \psi_\eps^N(x,\xi) \big)\mid_{\nu = \xi}
\eeq
where the notation $\F_{x \to \nu}(\ldots)\mid_{\nu = \xi}$ emphasizes that
the Fourier transform (in the $x$ variable) is carried out at fixed parameter
values $\xi$, $\eps$, and $N$, of its functional argument and then
evaluated at Fourier variable $\nu$ set equal to the parameter $\xi$.

Intuitively, rapid decrease of $J_\eps^N(\xi)$ would follow if we could
replace the family $\psi_\eps^N(.,\xi)$ by a single test function $\psi$ with
$\psi(x_0) = 1$. Note that for all $\eps\in (0,1]$, $\xi\in\Ga$,
$|\xi| \geq r_\eps$, and  $N\in\N$:
\beq\label{psi_supp}
    \supp(\psi_\eps^N(.,\xi)) \subseteq \supp(\vphi).
\eeq
The following condition specifies a regularity property of the family
$\psi_\eps^N(.,\xi)$ which will finally yield rapid decrease of
$J_\eps^N(\xi)$.
\begin{asmpt}\label{psi_asmpt}
 The family $\psi_\eps^N(.,\xi)$ ($\xi\in\Ga$, $|\xi| \geq r_\eps$,
    $N\in\N$,   $\eps\in (0,1]$)
    satisfies the following regularity condition: there is
    $M\in\N_0$, $0 \leq \de < 1$, and $\tau_0\in\R$ with the property that
    $\forall \al\in \N_0^n$ $\forall N\in\N$ $\exists C > 0, 0 < \eps_0 < 1$
    and $\exists r_\eps >0$ of slow scale such that
\beq \label{psi_reg}
    |\d_x^\al \psi_\eps^N(x,\xi)| \leq C \eps^{-M}\,
        (1+|\xi|)^{\de |\al|+\tau_0}
\eeq
for all $x\in\supp(\vphi)$, $\xi\in\Ga$ with $|\xi| \geq r_\eps$,
$\eps\in (0,\eps_0)$.
\end{asmpt}

\begin{lemma}\label{psi_lemma} Under Assumption \ref{psi_asmpt} let
    $\Ga_0\subset \Ga \cup \{ 0 \}$ be a closed conic neighborhood of $\xi_0$.
    Then $J_\eps^N(\xi)$ is rapidly decreasing when $\xi \in \Ga_0$.
\end{lemma}
\begin{proof} (This is a modified variant of a similar proof in \cite{HO:04},
the main difference being that in the present case the coupling of the
parameters $\xi$ and $\eps$ cannot be compensated for simply by homogeneity
arguments.)
Let $\chi\in\D(U)$ such that $\chi = 1$ on $\supp(\vphi) \supseteq
\supp(\psi_\eps^N)$ and let $\xi\in\Ga$, $|\xi| \geq r_\eps$. We have
\[
    |J_\eps^N(\xi)| = | \F_{x \to \nu}\big( (\chi f_\eps)(x) \,
        \psi_\eps^N(x,\xi) \big)\mid_{\nu = \xi} =
    \frac{1}{(2\pi)^n} |\big( g_\eps^N(.,\xi) * \FT{(\chi f_\eps)}\big)
        (\xi)|
\]
where $g_\eps^N(\eta,\xi) :=
\F(\psi_\eps^N(.,\xi))(\eta)$. This implies
\begin{multline*}
    (2\pi)^n |J_\eps^N(\xi)| \leq \int |g_\eps^N(\xi-\eta,\xi)|
        |\FT{(\chi f_\eps)}(\eta)| \, d\eta \\
    = \int_{\eta\in\Ga}  |g_\eps^N(\xi-\eta,\xi)|
        |\FT{(\chi f_\eps)}(\eta)| \, d\eta
        + \int_{\eta\in\Ga^c}  |g_\eps^N(\xi-\eta,\xi)|
        |\FT{(\chi f_\eps)}(\eta)| \, d\eta  \\
    =: J_{1,\eps}^N(\xi) + J_{2,\eps}^N(\xi).
\end{multline*}
Using Assumption \ref{psi_asmpt} we derive estimates on
$|g_\eps^N(\zeta,\xi)|$ as follows. For $\al\in\N_0^n$ arbitrary we have
\begin{multline*}
    |\zeta^\al g_\eps^N(\zeta,\xi)| =
        |\F( D_x^\al \psi_\eps^N(.,\xi))(\zeta)| \\
    \leq  \lone{D_x^\al \psi_\eps^N(.,\xi)}
    \leq C_\vphi \linf{D_x^\al \psi_\eps^N(.,\xi)}
    \leq C C_\vphi \eps^{-M} (1+|\xi|)^{\de|\al|+\tau_0}
\end{multline*}
where $M$ is independent of $\al$, $\eps < \eps_1$ as in (\ref{psi_reg}), and
$\zeta\in\R^n$ arbitrary.
Hence we have shown that for all $l\in\N_0$ $\exists C>0$ $\exists \eps_1 > 0$
such that
\beq\label{g_est}
    |g_\eps^N(\zeta,\xi)| \leq C \eps^{-M} (1 + |\zeta|)^{-l}
        (1+|\xi|)^{\de l + \tau_0}
\eeq
for all $\zeta\in\R^n$, $\eps < \eps_1$.

In estimating the integrand in $J_{1,\eps}^N(\xi)$ we use (\ref{g_est}),
the fact that $\FT{(\chi f_\eps)}$ is rapidly decreasing in $\Ga$, and
apply Peetre's inequality to obtain the following: there is $M_1$, $M_2$ such
that $\forall l,k \in \N_0$
\begin{multline*}
    |g_\eps^N(\xi-\eta,\xi)| |\FT{(\chi f_\eps)}(\eta)| \\
    \leq C_1 \eps^{-M_1} (1+|\xi|^2)^{(\tau_0+\de k)/2}
        (1 + |\xi-\eta|^2)^{-k/2}
        \eps^{-M_2} (1 + |\eta|^2)^{-l/2} \\
    \leq C' \eps^{-M_1-M_2} (1 + |\xi|^2)^{(\tau_0 - (1-\de) k)/2}
            (1 + |\eta|^2)^{(k-l)/2}
\end{multline*}
for suitable constants, $\eps$ sufficiently small, and $\eta\in\Ga$. If we
require $l > k + n$ then we may conclude that for arbitrary $k\in\N_0$
\[
    J_{1,\eps}^N(\xi) \leq C \eps^{-M_1-M_2} (1 + |\xi|)^{- (1-\de) k + \tau_0}
\]
where $M_1 + M_2$ is independent of $k$ and $\eps$ sufficiently small.

For a similar estimate of $J_{2,\eps}^N(\xi)$ we first note that $\FT{(\chi
f_\eps)}(\eta)$ is temperate in the following sense. There is $M_3$ and $C >
0$, $\eps_3 > 0$ such that
\[
    |\FT{(\chi f_\eps)}(\eta)| \leq C \eps^{-M_3} (1 + |\eta|^2)^{M_3/2}
\]
for all $\eta\in\R^n$ and $0 < \eps < \eps_3$. Applying this and (\ref{g_est}),
with $l+k$ instead of $l$, we obtain the following bound on the integrand in
$J_{2,\eps}^N(\xi)$:
\begin{multline*}
    |g_\eps^N(\xi,\xi-\eta)| |\FT{(\chi f_\eps)}(\eta)| \leq \\
     \leq C \eps^{-M_1-M_3} (1+|\xi|^2)^{(\tau_0 + \de(k+l))/2}
        (1 + |\xi - \eta|^2)^{-k/2}
        (1 + |\xi - \eta|^2)^{-l/2} (1 + |\eta|^2)^{M_3/2}.
\end{multline*}
By Peetre's inequality, $(1 + |\xi - \eta|^2)^{-k/2} \leq 2^k (1 + |\xi|^2)^{-k/2} (1
+ |\eta|^2)^{k/2}$. Furthermore, one can find $d> 0$ (resp.\ $d'>0$) such that
$\xi\in\Ga_0 $ and $\eta\in\Ga^c$ implies $|\xi - \eta| \geq d |\eta|$ (resp.\ $|\xi
- \eta| \geq d' |\xi|$) (cf.\ \cite[proof of Lemma 8.1.1]{Hoermander:V1-4}).
Therefore, we can write $(1+|\xi|^2)^{\de l/2} \leq C' (1 + |\xi-\eta|^2)^{\de l/2}$
and $(1 + |\xi - \eta|^2)^{-(1-\de)l/2} \leq  C (1 + |\eta|^2)^{-(1-\de)l/2}$,
showing that the integrand is bounded by some constant times
\[
    \eps^{-M_1-M_3} (1 + |\xi|)^{\tau_0 - (1- \de) k}
            (1 + |\eta|)^{k + M_3 - (1-\de) l}
\]
with $M_1$ and $M_3$ independent of $k$, $l$ and $\eps$ sufficiently small.
Requiring $(1-\de) l > k + M_3 + n$ yields
\[
    J_{2,\eps}^N(\xi) \leq C \eps^{-M_1-M_3} (1 + |\xi|)^{\tau_0 - (1-\de) k}.
\]
Hence we have proved rapid decrease of $J_\eps^N(\xi)$.
\end{proof}

To summarize the preceding discussion, Lemmas \ref{R_lemma} and
\ref{psi_lemma} imply the follow result.
\begin{prop}\label{basic_prop}
Let $P(x,D)$ be a partial differential operator with coefficients in
$\G(\Om)$, represented by the family $(P_\eps(x,D))_\eps$, and
$(R_\eps)_\eps$ be constructed according to Lemma \ref{alg_lem} (with
$A=P_\eps$ and $Q = \tP_\eps$). Let $u\in\G(\Om)$ and assume
$(x_0,\xi_0) \not\in \WF_g(Pu)$ with $U\times\Ga$, $\vphi$
as in (\ref{f_rap_dec}). Let  $(\psi_\eps^N)_{\eps,N}$ be defined by
(\ref{psi_def}). If $(R_\eps)_\eps$ and $(\psi_\eps^N)_{\eps,N}$ satisfy
Assumptions \ref{R_asmpt} and \ref{psi_asmpt} then $(x_0,\xi_0)\not\in
\WF_g(u)$.
\end{prop}

In the remainder of this paper we investigate various possibilities of
conditions on the operator family $(P_\eps)_\eps$, or its coefficients,
such that the crucial Assumptions \ref{R_asmpt} and \ref{psi_asmpt} are
guaranteed. In all these cases Proposition \ref{basic_prop} will allow us to
deduce microlocal regularity properties of a Colombeau solution to the equation
$P u = f$.

\section{Microlocal hypoellipticity conditions}

Throughout this section let $(P_\eps)_\eps$ be a family of linear partial
differential operators whose coefficients, $(a_\al^\eps)_\eps$, are
representatives of generalized functions in $\G(\Om)$. Denote by $P$ the
corresponding operator on $\G(\Om)$, mapping $\cl{(u_\eps)_\eps}$ into
$\cl{(P_\eps u_\eps)_\eps}$.

\begin{lemma} Let $P$ have coefficients in $\G^\infty(\Om)$. Then for any
$u\in\G(\Om)$
\beq
    \WF_g(P u) \subseteq \WF_g(u).
\eeq
\end{lemma}
\begin{proof}
$\WF_g(D^\al u) \subset \WF_g(u)$ is clear from the properties of the
Fourier transform. Furthermore, if $a\in\G^\infty$ then $\WF_g(a u)
\subseteq \WF_g(u)$ is a special case of \cite[Theorem 3.1]{HK:01}.
\end{proof}

\begin{thm}\label{mh_thm}
Let $P$ be a partial differential operator of order $m$ with coefficients
in $\G(\Om)$. Let $(x_0,\xi_0)\in \CO$ with an open conic neighborhood
$U\times\Ga$ and $m_0\in\R$, $0 \leq \de < \rho \leq 1$ such that the
following hypotheses are satisfied for any compact subset
$K \Subset U$:
\begin{enumerate}
\item $\exists q > 0$ $\exists r_\eps > 0$ of slow scale and $\exists \eps_0 >
    0$ such that
    \beq\label{mh_1}
        |P_\eps(x,\xi)| \geq \eps^q \, (1+|\xi|)^{m_0}
    \eeq
    for all $(x,\xi)\in K\times\Ga$, $|\xi| \geq r_\eps$, and $0<\eps<\eps_0$.
\item $\forall \al\in\N_0^n$ $\exists s_\eps^\al, r_\eps^\al > 0$ of slow
    scale and $\eps_\al > 0$ such that for all $\be\in\N_0^n$ with
    $0\leq|\be|\leq m$
    \beq\label{mh_2}
        |\d_x^\al \d_\xi^\be P_\eps(x,\xi)| \leq
            s_\eps^\al \, |P_\eps(x,\xi)|\, (1 + |\xi|)^{\de|\al|-\rho|\be|}
    \eeq
    for all $(x,\xi)\in K\times\Ga$, $|\xi| \geq r_\eps^\al$,
    and $0<\eps<\eps_\al$.
\end{enumerate}
Then we have, with $\La = U\times\Ga$, for any $u\in\G(\Om)$
\beq\label{ml_reg0}
\La \cap \WF_g( u) =  \La \cap \WF_g(P u).
\eeq
\end{thm}

\begin{rem}\label{mh_ginf_rem} Note that condition (\ref{mh_2}) implies that
all coefficients of $P$ are $\Ginf$ on $U$. Indeed, by the moderateness of
$P_\eps(x,\xi)$ and the fact that $s_\eps^\al \leq \eps^{-1}$ one obtains
uniform $\eps$-growth for all derivatives of $P_\eps$. Then one makes use of
the polynomial structure with respect to $\xi$ to first extract each highest
order coefficient separately (i.e., $|\be| = m$) and directly deduces its
regularity. Finally proceeding successively to lower orders, i.e., $|\be| = m
- 1$ and so on, each of the coefficients appears as the only one of the
current order with additional linear combinations of higher order
coefficients. Thus the regularity follows.
\end{rem}

The following statement is an immediate consequence, restating (\ref{ml_reg0}) as an inclusion relation.
\begin{cor}\label{mh_cor}
 Let $M_g(P)$ be the union of all open conic subsets $\La \subseteq
\CO$ where $P$ satisfies (\ref{mh_1}-\ref{mh_2}). Then the following
inclusion relation holds for any generalized function $u\in\G(\Om)$:
\beq \label{ml_reg}
    \WF_g(u) \subseteq \WF_g(P u) \cup M_g(P)^c.
\eeq
\end{cor}

\emph{Proof of Theorem \ref{mh_thm}.} We have to show that the families of
operators $R_\eps$ and functions $\psi_\eps^N(.,\xi)$  constructed in the
previous section satisfy Assumptions \ref{R_asmpt} and \ref{psi_asmpt}. Then the assertion follows from Proposition \ref{basic_prop}. Note
that (\ref{mh_1}) guarantees that $P_\eps(x,\xi)$ is staying away from zero in
the regions considered and hence the constructions according to Lemma
\ref{alg_lem} are well-defined.

Equations (\ref{r0}-\ref{rb}) define the coefficients, $r_0^\eps$ and
$r_\be^\eps$, of $R_\eps$ and show that we have to give appropriate bounds of
the generic factors which appear in the coefficients of powers of the
operator $R_\eps$:
\begin{align}
    & \d_x^\al\Big(\frac{\tP_\eps(x,-\xi)-P_\eps(x,\xi)}{P_\eps(x,\xi)}\Big)
        \label{gen_term_1} \\
    & \d_x^\al\Big( \d_\xi^{\ga + \be}\big(\tP_\eps(x,-\xi)\big)
    D_x^\ga\big(\frac{1}{P_\eps(x,\xi)}\big)\Big). \label{gen_term_2}
\end{align}

\emph{Step 1:} For each $K \Subset \Om$ and $\al,\be,\ga\in\N_0^n$
arbitrary $\exists S_\eps^{\al,\ga+\be}, p_\eps^{\al,\ga+\be} > 0$ of slow
scale and $\mu_{\al,\ga+\be}>0$ such that
\beq\label{step_1}
    |\d_x^\al\d_\xi^{\ga + \be}\big(\tP_\eps(x,-\xi)\big)|
    \leq S_\eps^{\al,\ga+\be}\,
        |P_\eps(x,\xi)| \, (1 + |\xi|)^{-\rho |\ga+\be| + \de|\al|}
\eeq
for all $(x,\xi)\in K\times\Ga$, $|\xi| \geq p_\eps^{\al,\ga+\be}$,
$0 < \eps < \mu_{\al,\ga+\be}$.

To see this we use the symbol expansion
\[
    \tP_\eps(x,\eta) = \sum_{0\leq|\sig|\leq m} \frac{(-1)^{|\sig|}}{\sig!}
        \big(\d_\xi^\sig\d_x^\sig P_\eps\big)(x,-\eta)
\]
which gives
\[
    \d_x^\al \d_\xi^{\ga + \be}\big(\tP_\eps(x,-\xi)\big) =
    \sum_{0\leq|\sig|\leq m-|\ga+\be|} \frac{(-1)^{|\sig|}}{\sig!}
        \big(\d_\xi^{\sig+\ga+\be} \d_x^{\sig+\al} P_\eps\big)(x,\xi).
\]
Applying (\ref{mh_2}) term by term and choosing a slow scale net
$S_\eps^{\al,\ga+\be}$ which dominates all appearing constants and slow scale
factors, as well as choosing $p_\eps^{\al,\ga+\be}$ to be the maximum of the
occurring radii 
$r_\eps^{\sig+\al}$, the assertion (\ref{step_1}) is immediate.

\emph{Step 2:} For each $K \Subset \Om$ and $\la\in\N_0^n$ arbitrary
$\exists T_\eps^{\la}, t_\eps^{\la} > 0$ of slow scale and $\nu_{\la}>0$ such
that
\beq\label{step_2}
    |D_x^\la\big(\frac{1}{P_\eps(x,\xi)}\big)| \leq
    T_\eps^{\la}\, |\frac{1}{P_\eps(x,\xi)}|\, (1+|\xi|)^{\de|\la|}
\eeq
for all $(x,\xi)\in K\times\Ga$, $|\xi| \geq t_\eps^{\la}$,
$0 < \eps < \nu_{\la}$.

The assertion is trivial if $|\la|=0$, so we assume $|\la|\geq 1$ and proceed
by induction. Differentiating the equality $ 1 = P_\eps/P_\eps$ we obtain,
by Leibniz' rule,
\[
    0 = P_\eps\cdot \d_x^\la(\frac{1}{P_\eps}) +
        \sum_{{\sig \leq \la}\atop {0\leq|\sig|<|\la|}}
        \binom{\la}{\sig}
        \d_x^{\la-\sig}P_\eps\cdot\d_x^\sig(\frac{1}{P_\eps}).
\]
Here, (\ref{step_2}) is applicable to each term in the sum over $\sig$ and
combination with (\ref{mh_2}) yields
\[
    |D_x^\la\big(\frac{1}{P_\eps(x,\xi)}\big)| \leq
    \sum_{{\sig \leq \la}\atop {0\leq|\sig|<|\la|}}
        \binom{\la}{\sig}
    s_\eps^{\la-\sig} T_\eps^{\sig}\, (1+|\xi|)^{\de(|\la-\sig|+|\sig|)}
\]
when $(x,\xi)\in K\times\Ga$, $|\xi| \geq \max_{\la-\sig}(r^{\la-\sig}_\eps)$,
and $\eps$ sufficiently small. From this we see that $T_\eps^{\la}$,
$t_\eps^{\la}$, and $\nu_{\la}$ can be chosen appropriately under a finite
number of conditions so that (\ref{step_2}) can be satisfied.

\emph{Step 3:} For each $K \Subset \Om$ and $\al\in\N_0^n$
$\exists c_\eps^{\al}, a_\eps^{\al} > 0$ of slow scale and $\mu_{\al}>0$ such
that for all  $1 \leq |\be| \leq m$
\beq\label{step_3}
    |\d_x^\al r_\be^\eps(x,\xi)| \leq c_\eps^{\al}\,
        (1+|\xi|)^{-\rho|\be| + \de|\al|}
\eeq
for all $(x,\xi)\in K\times\Ga$, $|\xi| \geq a_\eps^{\al}$,
$0 < \eps < \mu_{\al}$.

We have to estimate terms (\ref{gen_term_2}) when $\be\not=0$, which according
to Leibniz' rule are linear combinations of terms ($\sig \leq \al$)
\[
    \d_x^\sig \d_\xi^{\ga + \be}\big(\tP_\eps(x,-\xi)\big)\cdot
        D_x^{\al-\sig+\ga}\big(\frac{1}{P_\eps(x,\xi)}\big).
\]
Combining (\ref{step_1}) and (\ref{step_2}) gives upper bounds, for $|\xi|$
larger than some slow scale radius, of the form of some slow scale net times
$(1+|\xi|)^{-\rho |\ga+\be| + \de |\sig| + \de |\al-\sig+\ga|}$ which has
exponent $-(\rho - \de)|\ga| - \rho|\be| + \de |\al|$ and proves the
assertion since $\rho > \de$. (The appropriate slow scale nets are chosen, for
each $\al$, subject to finitely many conditions; and this may be done
uniformly over $1\leq |\be|\leq m$.)

\emph{Step 4:}  For each $K \Subset \Om$ and $\al\in\N_0^n$
$\exists d_\eps^{\al}, b_\eps^{\al} > 0$ of slow scale and $\nu_{\al}>0$ such
that
\beq\label{step_4}
    |\d_x^\al r_0^\eps(x,\xi)| \leq d_\eps^{\al}\,
        (1+|\xi|)^{-(\rho-\de) + \de|\al|}
\eeq
for all $(x,\xi)\in K\times\Ga$, $|\xi| \geq b_\eps^{\al}$,
$0 < \eps < \nu_{\al}$.

We have to find bounds on (\ref{gen_term_2}) when $\be = 0$ but $|\ga| \geq 1$
and on (\ref{gen_term_1}). The first case is done as in Step 3 and yields as
upper bound a slow scale net times $(1+|\xi|)^{-(\rho-\de)  + \de |\al|}$
since now $|\ga|\geq 1$. The term (\ref{gen_term_1}) is a linear combination
of terms ($\la \leq \al$)
\[
    \d_x^\la \big(\tP_\eps(x,-\xi)- P_\eps(x,\xi)\big)\cdot
        D_x^{\al-\la}\big(\frac{1}{P_\eps(x,\xi)}\big).
\]
Here, a bound on the second factor has the usual slow scale data coming with
$(1+|\xi|)^{\de |\al-\la|}/|P_\eps(x,\xi)|$ according to (\ref{step_2}). By
the symbol expansion of $\tP_\eps(x,-\xi)$ the factor on the left is seen to
be a linear combination of the following terms with $|\sig|
\geq 1$
\[
    \d_\xi^\sig \d_x^{\la+\sig} P_\eps(x,\xi).
\]
We obtain upper bounds with slow scale radii and factors times
$|P_\eps(x,\xi)| (1+|\xi|)^{-(\rho-\de) |\sig| + \de |\la|}$. Hence this
yields, apart from similar slow scale data, a common bound
$(1+|\xi|)^{-(\rho-\de) + \de |\la|}$ for the first factor in the product
above since $|\sig|\geq 1$. Multiplication of the bounds on both factors gives
finally the asserted upper bound. (We note once more that all required slow
scale nets can be chosen subject to finitely many conditions at fixed $\al$.)

\emph{Step 5:} $R_\eps(\xi;x,D)$ satisfies Assumption \ref{R_asmpt} with
$\tau = \rho - \de$ and $M_1 = 1$.

We prove this by induction. If $N=1$ it follows directly from (\ref{step_3})
and (\ref{step_4}) (set $\al=0$ in both equations) that, with some slow scale
net $s_\eps^0$,
\[
    |R_\eps(\xi;x,D)\vphi(x)| \leq s_\eps^0\, (1+|\xi|)^{-(\rho-\de)}
\]
when $(x,\xi)\in \supp(\vphi)\times\Ga$, $|\xi|$ larger than some slow scale
radius, and $\eps$ sufficiently small. Assume that for $j = 1,\ldots,N$ we
have, with some slow scale nets $s_\eps^j$, the induction hypothesis
\beq\label{ind_hyp}
    |R_\eps(\xi;x,D)^j \vphi(x)| \leq s_\eps^j\,
        (1+|\xi|)^{-j \tau}
\eeq
under similar conditions as above on $x$, $\xi$, and $\eps$. We let a term
$r_\be^\eps D_x^\be$ ($0\leq|\be|\leq m$) act on $R_\eps^N \vphi$ from the
left. Any derivative $D_x^\la$, $\la\leq \be$, falling on derivatives of
$r_\sig^\eps$ or $\vphi$ raises, according to (\ref{step_3})
and (\ref{step_4}), an overall upper bound at most by some slow
scale factor times $(1+|\xi|)^{\de |\be|}$. Again by (\ref{step_3})
and (\ref{step_4}), the additional factor $r_\be^\eps$
then brings in another slow scale net times $(1+|\xi|)^{-\rho |\be|}$,
if $\be \not= 0$, or $(1+|\xi|)^{-(\rho-\de)}$, if $\be = 0$. In any case,
the new $\xi$-exponents, added at this stage, sum up to at least $-(\rho-\de)$
and all slow scale factors and radii can be chosen subject to finitely many
conditions when $N+1$ is fixed. Combining this with the bounds on the terms in
$R_\eps^j \vphi$ we obtain estimate (\ref{ind_hyp}) with $N+1$ instead of $j$.
In particular, we observe that all appearing slow scale factors can be
compensated for by $1/\eps$, yielding the assertion.

\emph{Step 6:} $\psi_\eps^N(x,\xi)$ satisfies Assumption \ref{psi_asmpt} with
$M=q + 1$, $\de$ from (\ref{mh_2}), and $\tau_0 = -m_0 - \tau$.

The case $N=1$ is trivial, hence we assume $N \geq 2$.
Recalling (\ref{psi_def}) we rewrite $\d_x^\al(\psi_\eps^N(x,\xi))$ as linear
combination of the following terms ($\sig \leq \al$)
\[
    \d_x^\sig(w_\eps^N(x,\xi)) \cdot
        \d_x^{\al - \sig}(\frac{1}{P_\eps(x,\xi)}).
\]
Thanks to (\ref{step_2}) and (\ref{mh_1}) the second factor has a bound
$\eps^{-q} T_\eps^{\al-\sig} (1 + |\xi|)^{\de (|\al-\sig|) - m_0}$, on the
usual domains for $x$ and $\xi$ when $\eps$ is small.

Considering (\ref{w_def}) we can argue in a similar way as in Step 5 that
$\d_x^\sig$ acting on any term $R_\eps^k \vphi$ ($k=1,\ldots,N$) raises the
$\xi$-power in its overall upper bound at most by $\de|\sig|$. Together with
the bound (from (\ref{ind_hyp})) of the form slow scale times $(1+|\xi|)^{-\tau}$
of $\sum_{1 \leq k \leq N} |R_\eps^k \vphi|$ we obtain all in all the upper
bound, for some slow scale net $s_\eps^\al$,
\[
    |\d_x^\al(\psi_\eps^N(x,\xi))| \leq \eps^{-q} s^\al_\eps
        (1+|\xi|)^{-m_0 - (N-1) \tau + \de |\al|}
\]
when $x\in\supp(\vphi)$, $\xi\in\Ga$ with $|\xi|$ above some slow scale
radius, and $\eps$ sufficiently small. (Here, the choices of appropriate slow
scale nets are restricted by finitely many conditions at fixed $\al$.)
\ \hfill $\square$

\section{Some special cases and applications}

\paragraph{WH-Ellipticity with slow scales:} Slightly generalizing
a notion from \cite{HO:04}, an operator $P$ with Colombeau coefficients on
$\Om$ is said to be \emph{WH-elliptic with slow scales} if
on any compact subset $K \Subset \Om$ the following is valid:
\begin{enumerate}
\item $\exists q > 0$ $\exists r_\eps > 0$ of slow scale and $\exists \eps_0 >
    0$ such that
    \beq\label{wh_1}
        |P_\eps(x,\xi)| \geq \eps^q \, (1+|\xi|)^m
    \eeq
    for all $(x,\xi)\in K\times\R^n$, $|\xi| \geq r_\eps$, and $0<\eps<\eps_0$.
\item $\forall \al\in\N_0^n$ $\exists s_\eps^\al, r_\eps^\al > 0$ of slow
    scale and $\eps_\al > 0$ such that for all $\be\in\N_0^n$ with
    $0\leq|\be|\leq m$
    \beq\label{wh_2}
        |\d_x^\al \d_\xi^\be P_\eps(x,\xi)| \leq
            s_\eps^\al \, |P_\eps(x,\xi)|\, (1 + |\xi|)^{-|\be|}
    \eeq
    for all $(x,\xi)\in K\times\R^n$, $|\xi| \geq r_\eps^\al$,
    and $0<\eps<\eps_\al$.
\end{enumerate}
The following consequence of Theorem \ref{mh_thm} is an elliptic regularity
result (cf.\ \cite{HO:04} for related results in the constant coefficient
case).
\begin{cor} Let $P$ be WH-elliptic with slow scales. Then for all $u\in\G(\Om)$
\beq\label{ell_reg}
    \WF_g(u) = \WF_g(Pu).
\eeq
\end{cor}

\paragraph{First-order operators with slow scale coefficients:}
Here it is possible to obtain microlocal regularity from estimates of the
principal part over conic regions. Operators of this type were considered
earlier in \cite{HdH:01,HK:01}. Let $P_\eps(x,\xi) = \sum_{j=1}^n a_j^\eps(x)
\xi_j + a_0^\eps(x)$ with $a_k^\eps$ ($k=0,\ldots,n$) having slow scale
$\eps$-growth on compact sets in each derivative. Let $(x_0,\xi_0)\in\CO$ and
assume that $U\times\Ga$ is an open neighborhood where the following holds:
for each $K\Subset\Om$ there are $s_\eps, r_\eps > 0$ of slow scale and
$\eps_0 > 0$ such that
\beq\label{1st_order_est}
    |P_{\eps,1}(x,\xi)| \geq \frac{1}{s_\eps} (1+|\xi|)
\eeq
for all $(x,\xi)\in K\times\Ga$, $|\xi| \geq r_\eps$, and $0 < \eps < \eps_0$.
\begin{prop}\label{1st_order_prop}
Let $P$ be a first-order operator with variable slow scale Colombeau
coefficients and $(x_0,\xi_0)\in\CO$ be such that property
(\ref{1st_order_est}) holds. Then for any $u\in\G(\Om)$,
$(x_0,\xi_0)\not\in\WF_g(Pu)$ implies $(x_0,\xi_0)\not\in\WF_g(u)$.
\end{prop}
\begin{proof}
We show that $P$ satisfies (\ref{mh_1}-\ref{mh_2}) (with $m_0 = 1$, $\de = 0$, and
$\rho = 1$) when $(x,\xi)\in K\times\Ga$. First, using the slow scale property of
$a_0^\eps$ and (\ref{1st_order_est}) we obtain, with some slow scale net $s_\eps^0$,
\[
    |P_\eps(x,\xi)| \geq |P_{\eps,1}(x,\xi)| - |a_0^\eps(x)|
    \geq \frac{1+|\xi|}{s_\eps} - s_\eps^0 \geq \frac{1+|\xi|}{2 s_\eps}
\]
if $|\xi| \geq 2 s_\eps s_\eps^0$ and $\eps$ sufficiently small. This is
(\ref{mh_1}).

There are only first order nontrivial $\xi$-derivatives, and we have with some
slow scale net $s_\eps^j$
\[
    |\d_{\xi_j} P_\eps(x,\xi)| = |a_j^\eps(x)| \leq s_\eps^j =
        \frac{s_\eps^j 2 s_\eps (1+|\xi|)}{(1+|\xi|) 2 s_\eps} \leq
        2 s_\eps^j  s_\eps |P_\eps(x,\xi)| (1+|\xi|)^{-1}
\]
where we have used the above estimate on $|P_\eps(x,\xi)|$, and assume that $|\xi|$
is larger than some slow scale radius and $\eps$ small. Finally, the estimate of the
$x$-derivatives is also straight forward. Let $\al\in\N_0^n$ and assume that
$|\d_x^\al a_k^\eps(x)|\leq s_k^\eps$ (slow scale) when $x\in K$. Then we have
\begin{multline*}
    |\d_x^\al P_\eps(x,\xi)| \leq
        \sum_{j=1}^n |\d_x^\al a_j^\eps(x)| |\xi| + |\d_x^\al a_0^\eps(x)|\\
    \leq C (\max_{k=0,\ldots,n} s_\eps^k)\, (1+|\xi|)
        \leq s'_\eps\, |P_\eps(x,\xi)|
\end{multline*}
where $s'_\eps$ can be chosen to be $2C s_\eps \max s_\eps^k$ and the usual
assumptions on $|\xi|$ and $\eps$ are in effect.
\end{proof}

\paragraph{Conditions on the principal part:}

Let $P$ be an operator of order $m$ with Colombeau coefficients $a_\be^\eps(x)$. Let
$(x_0,\xi_0)\in\CO$ and with an open conic neighborhood $U\times\Ga$ on which the
following holds: For all $K\Subset U$
\begin{enumerate}
\item $\exists s_\eps, r_\eps > 0$ of slow scale $\exists \eps_0 >0$ such that
\beq\label{st_1}
    |P_{\eps,m}(x,\xi)| \geq
        \frac{1}{s_\eps}\, \sum_{|\al|=m} |a_\al^\eps(x)| \cdot (1+|\xi|)^m
\eeq
when $(x,\xi)\in K\times\Ga$, $|\xi| \geq r_\eps$, and $0 < \eps < \eps_0$
\item $\forall\ga\in\N_0^n$ $\exists s_\eps^\ga, r_\eps^\ga >0$ of slow scale
and $\exists \eps_\ga >0$ such that for all $\be\in\N_0^n$,
$0\leq |\be| \leq m$
\beq\label{st_2}
    |\d_x^\ga a_\be^\eps(x)| \leq
        s_\eps^\ga \, \sum_{|\al|=m} |a_\al^\eps(x)|
\eeq
when $(x,\xi)\in K\times\Ga$, $|\xi| \geq r_\eps^\ga$, and
$0 < \eps < \eps_\ga$.
\end{enumerate}

\begin{lemma}\label{st_mh_2}
Hypotheses (\ref{st_1}) and (\ref{st_2}) imply (\ref{mh_2}).
\end{lemma}
\begin{proof} Denote $b_m^\eps(x) = \sum_{|\al|=m} |a_\al^\eps(x)|$.  Since
\[
    \d_x^\ga \d_\xi^\be P_\eps(x,\xi) =
    \sum_{\be \leq \sig, |\sig|\leq m} \d_x^\ga a_\sig^\eps(x)
        \frac{\sig!}{(\sig-\be)!} \xi^{\sig-\be}
\]
we obtain
\[
    |\d_x^\ga \d_\xi^\be P_\eps(x,\xi)| \leq
        C |\xi|^{m-|\be|} \!\!
        \sum_{\be \leq \sig, |\sig|\leq m}\!\! |\d_x^\ga a_\sig^\eps(x)|
        \leq C'\, s_\eps^\ga \, b_m^\eps(x)\, (1+|\xi|)^{m-|\be|}.
\]
On the other hand,
\begin{multline}\label{st_inv}
    |P_\eps(x,\xi)| \geq |P_{\eps,m}(x,\xi)| -
            |\xi|^{m-1}\sum_{|\al| \leq m-1} |a_\al^\eps(x)|  \\
    \geq b_m^\eps(x) \, C \, (1+|\xi|)^m\,
        \big(\frac{1}{s_\eps} - \frac{s_\eps^0}{1+|\xi|} \big)
     \geq \frac{C}{2 s_\eps} \, (1+|\xi|)^m \, b_m^\eps(x)
\end{multline}
if $|\xi| \geq 2 s_\eps s_\eps^0 $. Combining the two estimates above yields
(\ref{mh_2}).
\end{proof}

\begin{lemma}\label{st_mh_1}
Assume that one of the following equivalent conditions holds:

(i)  for all $L\Subset \Om$ there is $p\in\N$ and $\eps_1 > 0$
such that
\beq\label{inv_asmpt}
    \inf_{x\in L} \sum_{|\al|=m} |a_\al^\eps(x)| \geq \eps^p
        \qquad 0 < \eps < \eps_1,
\eeq
that is, $\sum_{|\al|=m} |a_\al|$ is invertible in $\G(U)$;

(ii) $\sum_{|\al|=m} |a_\al|^2$ is invertible in $\G(U)$. 

Then hypotheses (\ref{st_1}) and (\ref{st_2}) imply (\ref{mh_1}) as well.
\end{lemma}
\begin{proof}
That (i) and (ii) are equivalent follows from the
equivalence of the $l^1$- and the $l^2$-norm on $\C^m$. By \cite[Theorem
1.2.5]{GKOS:01} there is $p\in\N$ and $\eps_1 >0$ such that $b_m^\eps(x) =
\sum_{|\al|=m} |a_\al^\eps(x)| \geq \eps^p$ for all $x\in L$, $0 < \eps < \eps_1$.
Now (\ref{mh_1}) follows directly from (\ref{st_inv}).
\end{proof}
We summarize the above results in the following statement.

\begin{prop}\label{st_prop}  Let $P$ be a partial differential operator with
Colombeau coefficients on $\Om$. Assume that $\sum_{|\al|=m} |a_\al|^2$
is invertible in $\G$ on $U$ and let $(x_0,\xi_0)\in\CO$ and
(\ref{st_1}-\ref{st_2}) be satisfied in some open
conic neighborhood $U\times\Ga$. Then $P$ satisfies hypotheses
(\ref{mh_1}-\ref{mh_2}) in the same region. In particular, we have the
microlocal regularity property (\ref{ml_reg0}) for each $u\in\G(\Om)$.
\end{prop}

\begin{rem} 
(i) The invertibility assumption in Proposition \ref{st_prop} cannot
    be dropped in general. For example, consider the zero divisor
    $c=\cl{(c_\eps)_\eps}\in\tilde{\C}$, defined by $c_\eps = 0$,
    if $1/\eps\in\N$, and $c_\eps = i$ otherwise. Then the
    operator with symbol $P_\eps(\xi) = c_\eps \xi$ satisfies
    (\ref{st_1}-\ref{st_2}) on all of $\R\times\R\zs$ but $P u =
    0$ admits any non-regular solution of the following kind:
    let $u_\eps$ equal a representative of some element in
    $\G(\R) \setminus \Ginf(\R)$ if $1/\eps\in\N$, and $u_\eps = 0$
    otherwise. Note that $P$ is not hypoelliptic and does not satisfy (\ref{mh_1}).

(ii) On the other hand, the invertibility of $\sum_{|\al|=m} |a_\al|^2$ is not
necessary for hypoellipticity of an operator. Neither is it necessary for (the
stronger) conditions (\ref{mh_1}-\ref{mh_2}) to hold. For example, consider
$P_\eps(\xi) = a_\eps \xi + b_\eps$, where $a_\eps = 0$ if $1/\eps\in\N$, $a_\eps =
1$ otherwise, and $b_\eps = 1 - a_\eps$ for all $\eps$. Then one easily verifies that
(\ref{mh_1}-\ref{mh_2}) hold (e.g., with $m_0 = 0$, $q = 0$, $\eps_0 = 0$ $r_\eps =
2$, $\delta = 0$, $\rho = 1$, $s_\eps = 3$), whereas the principal part coefficient
$\cl{(a_\eps)_\eps}$ is not invertible. Furthermore, condition (\ref{st_2}) fails to
hold for $P$ while (\ref{st_1}) is trivially satisfied.

(iii) The operator with symbol $P_\eps(\xi) = \eps \xi + i$ satisfies
    (\ref{mh_1}-\ref{mh_2}) but not (\ref{st_2}) (since $1 = |i|
    \not\leq s_\eps^1 \eps$). However, note that in this example the
    invertibility assumption on the principal part coefficient is met. (Again, the
estimate (\ref{st_1}) is trivial.)
\end{rem}

\begin{ex}
In this example we consider the situation of a hyperbolic operator with
discontinuous coefficients. Such operators arise e.g.\ in acoustic wave
propagation in a medium with irregularly changing properties. Let $\rho$ be
 the density, $c$ the sound speed of the medium. The pressure (perturbation)
$p$ solves the equation $P(p) = 0$ where
\beq \label{ex_equ}
    P =  \d_t^2  - c(x)^2 \rho(x) \d_x\big(\frac{1}{\rho(x)} \d_x \big).
\eeq
A typical assumption on the medium properties is that both, $\rho$ and
$c$, are (time independent and) measurable functions varying between strictly
positive bounds (but allow, e.g., for jump discontinuities).
To illustrate our theory, we interpret the coefficients $\rho$ and $c$ as
elements of $\G(\R)$ with representatives satisfying
\[
    0 < \ga_0 \leq c_\eps(x) \leq \ga_1, \qquad
    0 < r_0 \leq \rho_\eps(x) \leq r_1
\]
for $x\in\R$ and $\eps\in (0,1]$. In the setting of $\G(\R^2)$, the equation
$P(p)=0$ can be uniquely solved even with generalized functions as initial
data (see \cite{O:89}), and thus can model propagation of strong disturbances
even in media with highly complex structure.

In such circumstances the regions of regularity of the solution provide valuable
information. While global propagation of regularity for the constant coefficient case
of (\ref{ex_equ}) was dealt with in \cite{HO:04}, we are now able to address the
general case here. A representative of the operator (\ref{ex_equ}) with coefficients
in $\G(\R^2)$ as above is given by $ P_\eps(x,t,\d_x,\d_t) = \d_t^2 - c_\eps(x)^2
\rho_\eps(x) \d_x\big( \rho_\eps(x)^{-1} \d_x \big)$ with symbol
\[
    P_\eps(x,t,\xi,\tau) = -\tau^2 + c_\eps(x)^2 \xi^2 -
        i c_\eps(x)^2 \frac{\rho_\eps'(x)}{\rho_\eps(x)} \xi.
\]
We are going to identify regions of microhypoellipticity for $P$, i.e., the
set $M_g(P)$ introduced in Corollary \ref{mh_cor}. First, note that $U
\times\Ga \subseteq M_g(P)$ implies that $(\rho\otimes 1)\mid_U$ and
$(c\otimes 1)\mid_U$ are $\Ginf$ by Remark \ref{mh_ginf_rem}. Hence, denoting
by
$$
S_g(c,\rho) := \big(\singsupp_g(c) \cup \singsupp_g(\rho)\big)\times \R
$$
the union of the singular supports of $c\otimes 1$ and $\rho\otimes 1$, we
have
\[
    S_g(c,\rho) \times \R^2\setminus \{0\} \subseteq M_g(P)^c .
\]
Let $U$ be open in $\R^2$ such that $U \cap S_g(c,\rho) = \emptyset$ (i.e.,
the coefficients of $P$ are $\Ginf$ on $U$).
We observe that in the regions of $\Ginf$-regularity, the coefficients
actually satisfy the stronger property of having slow scale growth in each
derivative. This follows from the boundedness of $c$ and $\rho$ by
Corollary \ref{sc_ginf_cor}. Together with the (constant) positive lower
bound it
yields that $1 / \rho$ has the same properties there.
It follows that $P_\eps$, restricted to $U
\times \R^2$, is of the structure
\[
    P_\eps(x,t,\xi,\tau) = -\tau^2 + c_\eps(x)^2 \xi^2 -
        i b_\eps(x) \xi
\]
where $b = \cl{(b_\eps)_\eps}$ is real and of slow scale in each derivative.

Let $0 < \theta < \ga_0$ and define the open cone $\Ga_\theta$ in $\R^2
\setminus \{ 0\}$ by the conditions $|\tau| < (\ga_0 - \theta) |\xi|$ or
$|\tau| > (\ga_1 + \theta)|\xi|$. We will show that $P$ is microhypoelliptic
on $U\times \Ga_\theta$.

In fact, we can apply Proposition \ref{st_prop} with the principal part
$P_{\eps,2}(x,t,\xi,\tau) = -\tau^2 + c_\eps(x)^2 \xi^2$. Clearly $1 + c_\eps(x)^2
\geq 1 + \ga_0^2$ is invertible; furthermore, all estimates required in (\ref{st_2})
are then trivially satisfied due to the slow scale properties of the coefficients. It
remains to check (\ref{st_1}). The two conditions defining $\Ga_\theta$ yield
immediately that $|P_{\eps,2}(x,t,\xi,\tau)| \geq (\ga_0^2 - (\ga_0 - \theta)^2)
(\xi^2 + \tau^2 /(\ga_0-\theta)^2)/2$, resp.\ $|P_{\eps,2}(x,t,\xi,\tau)| \geq
((\ga_1 + \theta)^2 - \ga_1^2) (\tau^2 + \xi^2 (\ga_1+\theta)^2)/2$; therefore,
$|P_{\eps,2}(x,t,\xi,\tau)| \geq d_\theta (1 + \xi^2 + \tau^2)$ when $|(\xi,\tau)|
\geq r_\theta$, for suitable positive constants $d_\theta$ and $r_\theta$. On the
other hand, $1 + c_\eps^2(x) \leq 1 + \ga_1^2$ hence (\ref{st_1}) follows easily.

Since $\theta$ was arbitrary in the interval $(0,\ga_0)$ we obtain, letting $\theta
\to 0$, that
\[
    S_g(c,\rho)^c \times W^c \subseteq M_g(P)
\]
where $W := \{ (\xi,\tau)\in\R^2 \mid \ga_0 |\xi| \leq |\tau|
    \leq \ga_1 |\xi| \}$. Taking complements and summarizing we have shown
    that
\[
    S_g(c,\rho) \times \R^2\setminus \{0\} \subseteq M_g(P)^c
        \subseteq \big( S_g(c,\rho) \times \R^2\setminus \{0\} \big)
        \cup \big( S_g(c,\rho)^c \times W \big).
\]
\end{ex}

\bibliographystyle{abbrv}
\bibliography{gueMO}

\end{document}